\documentclass[12pt]{article}
\usepackage{graphicx, color}
\usepackage{amssymb}
\usepackage{amsmath}
\usepackage{latexsym}
\makeatletter
\newcommand{\C}{\mathbb C}
\newcommand{\R}{\mathbb R}
\newcommand{\Z}{\mathbb Z}
\newcommand{\N}{\mathbb N}

\newtheorem{thm}{Theorem}[section]
\newtheorem{lem}[thm]{Lemma}
\newtheorem{pro}[thm]{Proposition}
\newtheorem{dfn}[thm]{Definition}

\newtheorem{cor}[thm]{Corollary}
\newtheorem{rk}[thm]{Remark}

\begin{document}
\\

\begin{center}
{\large \bf Continuum-wise expansiveness for generic
diffeomorphisms}

\vspace{1.0cm}

 Manseob Lee\\

\vspace{0.1cm}

  {\em    Department of Mathematics,
       Mokwon University, \\Daejeon, 302-729, Korea.}\\

\vspace{0.5cm}

       {\em E-mail:
        lmsds@mokwon.ac.kr.}\\

\end{center}

\medskip
\begin{abstract}

Let $M$ be a closed smooth manifold and let $f:M\to M$ be a
diffeomorphism.  $C^1$-generically, a continuum-wise expansive
satisfies Axiom A without cycles. Let $M=\mathbb{T}^3$ and let
$f:\mathbb{T}^3\to\mathbb{T}^3$. There are a $C^1$ neighborhood
$\mathcal{U}(f)$ of $f\in\mathcal{RT}(\mathbb{T}^3)$ and a
residual set $\mathcal{R}\subset\mathcal{U}(f)$ such that for any
$g\in\mathcal{R}$, $g$ is not continuum-wise expansive, where
$\mathcal{RT}(\mathbb{T}^3)$ is the set of all robustly transitive
diffeomorphisms on $\mathbb{T}^3.$

 \end{abstract}
\smallskip
\section{Introduction}
 Let $M$ be a closed smooth manifold with ${\rm
dim} M\geq2$, and let ${\rm Diff}(M)$ be the space of
diffeomorphisms of $M$ endowed with the $C^1$ topology. Denote by
$d$ the distance on $M$ induced from a Riemannian metric
$\|\cdot\|$ on the tangent bundle $TM$. In differentiable
dynamical systems,
  expansiveness is a very useful notion  to investigate for stability theory.
For instance, Ma\~n\'e \cite{M1} proved that   the $C^1$-interior
of the set of expansive diffeomorphisms coincides with the set of
quasi-Anosov diffeomorphisms. Here $f$ is {\it quasi-Anosov} if
for all $v\in TM\setminus\{0\}$, the set $\{\|Df^n(v)\|:n\in\Z\}$
is unbounded. Let $f\in{\rm Diff}(M)$.
 We say that $f$ is {\it expansive} if there is $e>0$ such that for
any $x, y\in M$ if $d(f^i(x), f^i(y))<e$ for all $i\in\Z$ then
$x=y.$ Denote by $\mathcal{E}$ the set of all expansive
diffeomorphisms. From now, we introduce various expansiveness
(N-expansive, countably expansive, measure expansive \cite{Mo,
MS}) which are general notions of original expansiveness.

 We say that $f$
is {\it N-expansive} if there is $e>0$ such that for any $x\in M$,
the number of elements of the set $\Gamma_{e}(x)=\{y\in M:
d(f^i(x), f^i(y))<e$ for all $i\in\Z\}$ is less than $N.$ Denote
by $\mathcal{GE}$ the set of all N-expansive diffeomorphisms on
$M$.
 We say that $f$ is {\it countably expansive} if there is
$e>0$ such that for $x\in M$, the number of elements of the set
$\Gamma_{e}(x)=\{y\in M: d(f^i(x), f^i(y))<e$ for all $i\in\Z\}$
is  countable, where $e$ is an expansive constant for $f.$

Note that if a diffeomorphism $f$ is expansive then
$\Gamma_{e}(x)=\{x\}$ for $x\in M.$ Thus if a diffeomorphism $f$
is  expansive then $f$ is countably expansive, but the converse is
not true (see \cite{MS}).

For a Borel probability measure $\mu$ on $M$, we say
that $f$ is $\mu$-expansive if there is $\delta>0$ such that
 $\mu(\Gamma_{e}(x))=0$ for all $x\in M.$ In this case, we say
that $\mu$ is a {\it  expansive measure} for $f$. We say that $f$
is {\it measure expansive} if it is $\mu$-expansive for every
non-atomic Borel probability measure $\mu$. Denote by
$\mathcal{ME}$ the set of all measure-expansive diffeomorphisms on
$M$.

Continuum-wise expansive diffeomorphisms was introduced by Kato
\cite{K}. A set $\Lambda$ is {\it nondegenerate} if the set
$\Lambda$ is not reduced to one point. We say that $\Lambda\subset
M$ is a {\it subcontinuum} if it is a compact connected
nondegenerate subset of $M$. \begin{dfn}  A diffeomorphism $f$ on
$M$ is said to be {\it continuum-wise expansive}  if there is a
constant $e>0$ such that for any nondegenerate subcontinuum
$\Lambda$ there is an integer $n=n(\Lambda)$ such that ${\rm diam}
f^n(\Lambda)\geq e$, where ${\rm diam} \Lambda = \sup\{d(x,y) : x,
y \in \Lambda \}$ for any subset $\Lambda\subset M$. Such a
constant $e$ is called a {\it continuum-wise expansive constant}
for $f$. \end{dfn} Note that every expansive diffeomorphism is
continuum-wise expansive diffeomorphism, but its converse is not
true (see \cite [Example 3.5]{K}). Denote by $\mathcal{CWE}$ the
set of all continuum-wise expansive diffeomorphisms of $M$.  In
\cite{A}, Artigue showed that
$$\mathcal{E}\Rightarrow \mathcal{GE}\Rightarrow
\mathcal{CE}=\mathcal{ME}\Rightarrow \mathcal{CWE},$$ where
$\mathcal{CE}$ is the set of all countably expansive
diffeomorphisms on $M$. For a $C^1$ perturbation expansive
diffeomorphism, we can find the following result (see \cite{ACO,
L, S1, SSY1}). Denote by $int A$ the $C^1$-interior of a set $A$
of $C^1$-diffeomorphisms of $M.$
\begin{thm} Let $f\in{\rm Diff}(M)$. Then we have the following
 $$int \mathcal{E}=int \mathcal{GE}= int\mathcal{ME}=int
\mathcal{CWE}.$$ \end{thm}

Let $\Lambda$ be a closed $f$-invariant set. We say that $\Lambda$
is {\it hyperbolic} if the tangent bundle $T_{\Lambda}M$ has a
$Df$-invariant splitting $E^s\oplus E^u$ and there exists
constants  $C>0$ and $0<\lambda<1$ such that
$$\|D_xf^n|_{E_x^s}\|\leq C\lambda^n\;\;{\rm and}\;\;\|D_xf^{-n}|_{E_x^u}\|\leq C\lambda^{n} $$
for all $x\in \Lambda$ and $n\geq 0.$ If $\Lambda=M$ then $f$ is
said to be Anosov.

 It is well know that
if a diffeomorphism $f$ is Anosov then it is quasi-Anosov, but the
converse is not true (see \cite{FR}). Thus if a diffeomorphism $f$
is Anosov then $f$ is expansive, N-expansive, measure expansive,
countably expansive and continuum-wise expansive. We say that $f$
satisfies {\it Axiom A} if the non-wandering set $\Omega(f)$ is
hyperbolic and it is the closure of $P(f).$ A point $x\in M$ is
said to be {\it non-wandering } for $f$ if for any non-empty open
set $U$ of $x$ there is $n\geq0$ such that $f^n(U)\cap
U\not=\emptyset.$ Denote by $\Omega(f)$ the set of all
non-wandering points of $f.$ It is clear that
$\overline{P(f)}\subset\Omega(f).$ A diffeomorphism $f$ is {\it
$\Omega$-stable} if there is a $C^1$-neighborhood $\mathcal{U}(f)$
of $f$ such that for any $g\in\mathcal{U}(f)$ there is a
homeomorphism $h:\Omega(f)\to\Omega(g)$ such that $h\circ f=g\circ
h,$ where $\Omega(g)$ is the non-wandering set of $g.$ A subset
$\mathcal{G}\subset {\rm Diff}(M)$ is called {\it residual} if it
contains a countable intersection of open and dense subsets of
${\rm Diff}(M)$. A dynamic property is called $C^1${\it generic}
if it holds in a residual subset of ${\rm Diff}(M).$ Arbieto
\cite{Ar} proved that  if a $C^1$ generic diffeomorphism $f$ is
expansive then it is Axiom A without cycles. Lee \cite{L} proved
that  if a $C^1$ generic diffeomorphism $f$ is N-expansive then it
is Axiom A without cycles. Very recently, Lee \cite{L2} proved
that if a $C^1$ generic diffeomorphism $f$ is measure expansive
then it is Axiom A without cycles. From that, we consider $C^1$
generic continuum-wise expansive diffeomorphisms.
The following is a main result.\\

 \noindent{\bf Theorem A} {\em For $C^1$ generic $f\in{\rm
Diff}(M),$ if $f$ is continuum-wise expansive then it is Axiom A
without
cycles.}\\


  We say that a $f$-invariant closed set $\Lambda$ admits a
{\it dominated splitting} if the tangent bundle $T_{\Lambda}M$ has
a continuous $Df$-invariant splitting $E\oplus F$ and there exist
constants $C>0$ and $0<\lambda<1$ such that
$$\|D_xf^n|_{E(x)}\|\cdot\|D_xf^{-n}|_{F(f^n (x))}\|\leq C\lambda^{n} $$
for all $x\in \Lambda$ and $n\geq 0.$ If the dominated splitting
can be written as a sum $$T_{\Lambda}M=E_1\oplus E_2\oplus\cdots
\oplus E_i\oplus E_{i+1}\oplus\cdots\oplus E_k,$$ then we say that
the sum is {\it dominated} if for all $i$ the sum $$(E_1\oplus
E_2\oplus\cdots\oplus E_i)\oplus (E_{i+1}\oplus E_{i+2}\oplus
\cdots E_k)$$ is dominated. Note that the decomposition is called
the {\it finest} dominated splitting if we can't decompose in a
non-trivial way subbundle $E_i$ appearing in the splitting.

The set $\Lambda$ is {\it partially hyperbolic} if there is a
dominated splitting $E\oplus F$ of $T_{\Lambda}M$ such that either
$E$ is contracting or $F$ is expanding.

\begin{dfn} We say that a compact $f$-invariant set
$\Lambda\subset M$ is strongly partially hyperbolic  if the
tangent bundle $T_{\Lambda}M$ has a dominated splitting $E^s\oplus
E^c\oplus E^u$ and there exist $C>0$ and $0<\lambda<1$ such that
for all $v\in E^s$, we have $\|Df^n(v)\|\leq C\lambda^n\|v\|$ for
all $n\geq0,$ and for all $v\in E^u$, we have $\|Df^{-n}(v)\|\leq
C\lambda^n\|v\|$ for all $n\geq0,$ where $E^c$ is the central
direction of the splitting.
\end{dfn}

Note that if $\Lambda$ is hyperbolic for $f$ then it is strongly
partially hyperbolic and $E^c$ is not empty, that is, $E^c=\{0\}.$
For a partially hyperbolic diffeomorphism, Burns and Wilkinson
\cite{BW} showed the following lemma.

\begin{lem}\label{fake} Let $\Lambda$ be a compact $f$-invariant
set with a partially hyperbolic splitting,
$$T_{\Lambda}M=E^s\oplus E_1^c\oplus\cdots\oplus E^c_k\oplus
E^u.$$ Let $E^{cs,i}=E^s\oplus E_1^c\oplus\cdots\oplus E_i^c$ and
$E^{cu, i}=E^c_i\oplus\cdots\oplus E^c_k\oplus E^u$ and consider
their extensions $\widetilde{E}^{cs, i}$ and
$\widetilde{E}^{cu,i}$ to a small neighborhood of $\Lambda.$ Then
for any $\epsilon>0$ there exist constants $R>r>r_1>0$ such that
for any $x\in\Lambda$, the neighborhood $B(x, r)$ is foliated by
foliations $\widehat{W}^u(x), \widehat{W}^s(x), \widehat{W}^{cs,
i}(x)$ and $\widehat{W}^{cu, i}(x)(i=1, \ldots, k)$ such that for
each $\sigma\in\{u, s, (cs, i), (cu, i)\}$ the following
properties hold.
\begin{itemize}
\item[(a)] {\em Almost tangency of invariant distributions.} For each
$y\in B(x, r)$, the leaf $\widehat{W}^{\sigma}_x(y)$ is $C^1$, and
the tangent space $T_y\widehat{W}^{\sigma}_x(y)$ lies in a cone of
radius $\epsilon$ about $\widetilde{E}^{\sigma}(y).$
\item[(b)] {\em Coherence.} $\widehat{W}^s_x$ subfoliates
$\widehat{W}^{cs, i}_x$ and $\widehat{W}^u_x$ subfoliates
$\widehat{W}^{cu, i}_x$ for each $i\in\{1, \ldots, k\}.$
\item[(c)] {\em Local invariance.} For each $y\in B(x, r)$ we have
$f(\widehat{W}^{\sigma}_x(y, r_1))\subset
\widehat{W}^{\sigma}_{f(x)}(f(y))$ and
$f^{-1}(\widehat{W}^{\sigma}_x(y, r_1))\subset
\widehat{W}^{\sigma}_{f^{-1}(x)}(f^{-1}(y)),$ where
$\widehat{W}^{\sigma}_x(y, r_1)$ is the connected components of
$\widehat{W}^{\sigma}_x(y)\cap B(y, r_1)$ containing $y.$
\item[(d)]{\em Uniqueness.} $\widehat{W}^s_x(x)=W^s(x, r)$ and
$\widehat{W}^u_x(x)=W^u(x, r).$
\end{itemize}
\end{lem}

We say that a diffeomorphism $f$ has a {\it homoclinic tangency}
if there is a hyperbolic periodic point $p$ whose invariant
manifolds $W^s(p)$ and $W^u(p)$ have a non-transverse
intersection. The set of $C^1$ diffeomorphisms that have some
homoclinic tangencies will be denoted $\mathcal{HT}.$ For a
homoclinic tangency, Pacifico and Vieitez \cite{PV} proved that
surface diffeomorphisms presenting homoclinic tangencies can be
$C^1$-approximated by non-measure expansive diffeomorphisms. Form
the result, Lee \cite{L1} proved that if $f$ has a homoclinic
tangency associated to a hyperbolic periodic point $p$, then there
is a $g$ $C^1$-close to $f$ such that $g$ is not continuum-wise
expansive.
\begin{pro}{\em \cite[Theorem 1.1]{CSY}} The diffeomorphism $f$ in
a dense $\mathcal{G}_{\delta}$ subset $\mathcal{G}\subset {\rm
Diff}(M)\setminus\overline{\mathcal{HT}}$ has the following
properties.

\begin{itemize}
\item[(a)] Any aperiodic class $\mathcal{C}$ is partially
hyperbolic with a one-dimensional central bundle. Morevoer, the
Lyapunov exponent along $E^c$ of any invariant measure supported
on $\mathcal{C}$ is zero.
\item[(b)] Any homoclinic class $H_f(p)$ has a partially hyperbolic
structure $$T_{H_f(p)}M=E^s\oplus E^c_1\oplus \cdots\oplus
E^c_k\oplus E^u.$$ Moreover, the minimal stable dimension of the
periodic orbits of $H_f(p)$ is ${\rm dim}E^s$ or ${\rm dim}E^s
+1.$ Similarly, the maximal stable dimension of the periodic
orbits of $H_f(p)$ is ${\rm dim}E^s+k$ or ${\rm dim}E^s+k-1.$ For
every $i=1, \ldots, k,$ there exist periodic points in $H_f(p)$
whose Lyapunov exponent along $E^c_i$ is arbitrarily close to 0.
In particular, if $f\in\mathcal{G}$, then $f$ is partially
hyperbolic.
\end{itemize}

\end{pro}
 Recently, Pacifico and Vieitez \cite{PV} proved that there is a
  residual subset $\mathcal{G}$ of ${\rm
Diff}(M)\setminus\overline{\mathcal{HT}}$  such that for any Borel
probability measure $\mu$ absolutely continuous  with respect to
Lebesgue, $f$ is $\mu$-expansive.
 Lee \cite{L2} showed that there is a partially hyperbolic
 diffeomorphism such that it is not measure expansive. From the
 facts, we consider continuum-wise expansive for partially
 hyperbolic diffeomorphisms. The set $\lambda$ is {\it
 transitive} if there is a point $x\in\Lambda$ such that
 $\omega(x)=\Lambda$, where $\omega(x)$ the omega limit set of
 $f.$ we say that the set $\Lambda$ is {\it robustly transitive }
 if there are a $C^1$-neighborhood $\mathcal{U}(f)$ of $f$ and a
 neighborhood $U$ of $\Lambda$ such that for any
 $g\in\mathcal{U}(f)$, $\Lambda_g(U)=\bigcap_{n\in\Z}g^n(U)$ is
 transitive, where $\Lambda_g(U)$ is the continuation of
 $\Lambda.$
 Let
$M=\mathbb{T}^3$ and let $f:\mathbb{T}^3
\to\mathbb{T}^3$ be a diffeomorphism.\\

\noindent{\bf Theorem B} {\em  There is a $C^1$ neighborhood
$\mathcal{U}(f)$ of $f\in{\mathcal{RT}}(\mathbb{T}^3)$ and a
residual set $\mathcal{R}\subset\mathcal{U}(f)$ such that for any
$g\in\mathcal{R}$,
 $g$ is not continuum-wise expansive, where $\mathcal{RT}(\mathbb{T}^3)$ is the set of all robustly transitive
diffeomorphisms on $\mathbb{T}^3$}.

\section{Proof of Theorems}
\subsection{Proof of Theorem A}

Let $M$ be as before, and let $f\in{\rm Diff}(M).$ The following
Franks' lemma \cite{F} will play essential roles
 in our proofs.

 \begin{lem}\label{frank} Let $\mathcal{U}(f)$
be any given $C^1$ neighborhood of $f$. Then there exist
$\epsilon>0$ and a $C^1$ neighborhood
$\mathcal{U}_0(f)\subset\mathcal{U}(f)$ of $f$ such that for given
$g\in \mathcal{U}_0(f)$, a finite set $\{x_1, x_2, \cdots, x_N\}$,
a neighborhood $U$ of $\{x_1, x_2, \cdots, x_N\}$ and linear maps
$L_i : T_{x_i}M\rightarrow T_{g(x_i)}M$ satisfying
$\|L_i-D_{x_i}g\|\leq\epsilon$ for all $1\leq i\leq N$, there
exists $\widehat{g}\in \mathcal{U}(f)$ such that
$\widehat{g}(x)=g(x)$ if $x\in\{x_1, x_2, \cdots,
x_N\}\cup(M\setminus U)$ and $D_{x_i}\widehat{g}=L_i$ for all
$1\leq i\leq N$.
\end{lem}

 The following
was proved by \cite{L2}. For convenience, we give the proof in the
section.

 \begin{lem}\label{lem21}
If $f\in {\rm Diff}(M)$ has a non-hyperbolic periodic point, then
for any neighborhood $\mathcal{U}(f)$ of $f$ and any $\eta>0$,
there are $g\in\mathcal{U}(f)$ and a curve $\gamma$ with the
following property:
\begin{enumerate}
\item $\gamma$ is $g$ periodic, that is, there is $n\in\mathbb{Z}$
such that $g^n(\gamma)=\gamma$;
\item the length of $g^i(\gamma)$ is less
than $\eta$, for all $i\in\mathbb{Z}$;
\item $\gamma$ is normally hyperbolic with respect to $g$.
\end{enumerate}
\end{lem}

\noindent{\bf Proof.} Let $\mathcal{U}(f)$ be a $C^1$ neighborhood
of $f$. Suppose $p$ is not hyperbolic periodic point of $f.$ For
simplicity, we may assume that $p$ is a fixed point of $f$. By
Lemma \ref{frank}, there is $g\in\mathcal{U}(f)$
 such that $D_pg^{\pi(p)}$ has an
eigenvalue $\lambda$ with $|\lambda|=1.$ Then $g(p)=p_g$ and
$T_{p_g}M=E^s_{p_g}\oplus E^c_{p_g}\oplus E^u_{p_g}$, where
$E^s_{p_g}$ is the eigenspace corresponding to the eigenvalues
with modulus less than 1,  $E^u_{p_g}$ is the eigenspace
corresponding to the eigenvalues with modulus more than 1, and
$E^c_{p_g}$ is the eigenspace corresponding to $\lambda$. If
$\lambda\in\R$ then ${\rm dim}E^c_{p_g}=1$ and if $\lambda\in \C$
then ${\rm dim}E^c_p=2.$

In the proof, we consider ${\rm dim}E^c_{p_g}=1.$ For case ${\rm
dim}E^c_{p_g}=2$, we can obtain the result as in the case ${\rm
dim}E^c_{p_g}=1$.

 Since ${\rm dim}E^c_{p_g}=1$ we assume that $\lambda=1$.
By Lemma \ref{frank}, there are $\epsilon>0$ and $h\in
\mathcal{U}(f)$ such that
\begin{itemize}
\item[$\cdot$] $h(p_g)=g(p_g)=p_g,$
\item[$\cdot$] $h(x)={\rm exp}_{p_g}\circ D_{p_g}g\circ {\rm
exp}_{p_g}^{-1}(x)$ if $x\in B_{\epsilon}(p_g),$ and
\item[$\cdot$] $h(x)=g(x)$ if $x\notin B_{4\epsilon}(p_g).$
\end{itemize}
Since $\lambda=1$, we can construct a closed  small arc
$\mathcal{I}_{p_g}\subset B_{\epsilon}(p_g)\cap {\rm
exp}_{p_g}(E^c_{p_g}(\epsilon))$ with its center at $p_g$ such
that
\begin{itemize}
\item[$\cdot$] ${\rm diam}\mathcal{I}_{p_g}=\epsilon/4$,
\item[$\cdot$] $h(\mathcal{I}_{p_g})=\mathcal{I}_{p_g},$ and
\item[$\cdot$] $h|_{\mathcal{I}_{p_g}}$ is the identity map.
\end{itemize}
Here $E^c_{p_g}(\epsilon)$ is the $\epsilon$-ball in $E^c_{p_g}$
centered at the origin $\overrightarrow{o}_p.$  Then
$\mathcal{I}_{p_g}$ is normally hyperbolic with respect to $h$,
and for any $\eta<\epsilon/4$, the length of $\mathcal{I}_{p_g}$
is less than $\eta$.
 \hfill$\square$\\

By the persistency of normally hyperbolic, we know that there is a
neighborhood $\mathcal{U}(g)$ of $g$ such that for any
$\tilde{g}\in\mathcal{U}(g)$ there is a curve $\tilde{\gamma}$
close to $\gamma$ such that all properties of $\gamma$ listed in
the Lemma \ref{lem21} is also satisfied for $\tilde{\gamma}$.


 For
$f\in{\rm Diff}(M)$, we say that $f$ is the {\it star
diffeomorphism} (or $f$ satisfies the {\it star condition}) if
there is a $C^1$-neighborhood $\mathcal{U}(f)$ of $f$ such that
all periodic points of $g\in\mathcal{U}(f)$ are hyperbolic. Denote
by $\mathcal{F}(M)$ the set of all star diffeomorphisms. Aoki
\cite{Ao} and Hayashi \cite{H} showed that   for any dimension
case, if $f\in\mathcal{F}(M)$ then $f$ is Axiom A without cycles.

\begin{lem}{\em \cite[Lemma 3.4]{L2}} \label{lem25}
There is a residual set $\mathcal{G}\subset{\rm Diff}(M)$ such
that for any $f\in\mathcal{G}$,
\begin{itemize}
\item [$\cdot$] either $f$ is a star,
\item [$\cdot$] or for any $\varepsilon>0$ there is a periodic curve $\gamma$ such that the length of
$f^n(\gamma)$ is less than $\varepsilon$, for any
$n\in\mathbb{Z}$.
\end{itemize}
\end{lem}

\noindent{\bf Proof.} Let $\mathcal{H}_n$ be the set of $C^1$
diffeomorphisms $f$ such that $f$ has a normally hyperbolic
$\gamma$ which is $1/n$-simply periodic curve. Since $\gamma$ is
normally hyperbolic, we know that $\mathcal{H}_n$ is open. Let
$\mathcal{N}_n={\rm Diff}(M)\setminus\overline{\mathcal{H}_n}.$
Then $\mathcal{H}_n(\eta)\cup\mathcal{N}_n(\eta)$ is open and
dense in ${\rm Diff}(M).$ Let
$\mathcal{G}=\bigcap_{n\in\N^+}(\mathcal{H}_n\cup\mathcal{N}_n).$
Then $\mathcal{G}$ is $C^1$ residual in ${\rm Diff}(M).$ Let $f\in
\mathcal{G}$ and assume $f$ is not a star diffeomorphism, we know
that $f\in\overline{\mathcal{H}_n}$ for any $n\in\N^+$ by Lemma
\ref{lem21}. Hence $f\notin \mathcal{N}_n$ and $f\in\mathcal{H}_n$
for any $n$. We know that $f$ has a normally hyperbolic $\gamma$
which is $\varepsilon$-simply periodic curve, for any
$\varepsilon>0$. \hfill$\square$\\

The following was proved by \cite[Lemma 2.2]{ACO}.

\begin{lem}\label{cw-exp}
  Let $C\subset M$ be a continuum. $f$ is continuum-wise
expansive if and only if there is $\delta>0$ such that for all
$x\in M$, if a continuum $C\subset \Gamma_{\delta}(x)$ then $C$ is
a singleton.
\end{lem}

\noindent{\bf Proof of Theorem A.} Let $f\in\mathcal{G}$ be
continuum-wise expansive. Suppose by contradiction that
$f\not\in\mathcal{F}(M).$ From Lemma \ref{lem25}, for any
$\varepsilon>0$ there is a periodic curve $\gamma$ such that the
length of $f^i(\gamma)$ is less than $\varepsilon$, for any
$i\in\mathbb{Z}$.   Let $\Gamma_{\epsilon}(x)=\{x\in M :d(f^i(x),
f^i(y))\leq \epsilon$ for all $i\in\Z\}.$ Since
$f^n(\gamma)=\gamma$ for some $n\in\Z$, we know $\gamma\subset
\Gamma_{\epsilon}(x).$ By Lemma \ref{cw-exp}, $\gamma$ should be a
singleton which is a contradiction since $\gamma$ is a nontrivial
continuum. Thus if $f\in\mathcal{G}$ is continuum-wise expansive
then it is Axiom A without cycles. \hfill$\square$\\

\subsection{Proof of Theorem B}

In this section, let $M=\mathbb{T}^3$ and let $f:\mathbb{T}^3
\to\mathbb{T}^3$ be a diffeomorphism. In \cite[Theorem B]{M2},
Ma\~n\'e constructed a robustly nonhyperbolic transitive
diffeomorphism $f\in{\rm Diff}(\mathbb{T}^3)$. By \cite[Theorem
B]{DPU}, every robustly transitive diffeomorphism $f$ on
$\mathbb{T}^3$ is partially hyperbolic. Thus we can find a
partially hyperbolic diffeomorphism $f$ on
$\mathbb{T}^3$ such that $f$ is robustly nonhyperbolic transitive. \\

 \begin{rk}
There is a partially hyperbolic diffeomorphism $f$ on
$\mathbb{T}^3$ such that $f$ is robustly nonhyperbolic
transitive.\end{rk}

\begin{lem} {\em \cite[Corollary D]{DPU} and \cite[Theorem
4.10]{ABC}} There is a residual set $\mathcal{G}_1\subset
\mathcal{RT}(\mathbb{T}^3)$ such that for any $f\in\mathcal{G}_1$,
$f$ is strongly partially hyperbolic, and $\mathbb{T}^3$ is the
homoclinic class $H_f(p)$, for some hyperbolic periodic point $p$.
\end{lem}
Let $M$ be a closed smooth $n(\geq2)$-dimensional manifold, and
let $f:M\to M$ be a diffeomorphism. The following notion was
introduced by Yang and Gan \cite{YG}. For any $\epsilon>0$, a
$C^1$ curve $\eta$ is called  a {\it $\epsilon$-simply periodic
curve} of $f$ if
\begin{itemize}\item[(i)]
 $\eta$ is diffeomorphic to $[0, 1]$ and its two
endpoints are hyperbolic periodic points of $f,$ \item[(ii)]
 $\eta$ is periodic with period $\pi(\eta)$ and
 $L(f^i(\eta))<\epsilon$ for any
$i\in\{1, 2, \cdots, \pi(\eta)\},$ where $L(\eta)$ denotes the
length of $\eta,$ and \item[(iii)] $\eta$ is normally hyperbolic.
\end{itemize}

\begin{lem}{\em \cite[Lemma 2.1]{YG}}\label{curve} There is a residual set $\mathcal{G}_2\subset{\rm Diff}(M)$ such that for any $f\in\mathcal{G}_2$, and any hyperbolic
periodic point $p$ of $f,$ we have the following:

 for any $\epsilon>0$, if for any $C^1$ neighborhood
$\mathcal{U}(f)$ of $f$ some $g\in\mathcal{U}(f)$  has a
$\epsilon$-simply periodic curve $\eta$ such that two endpoints of
$\eta$ are homoclinically related with $p_g$ then $f$ has an
$2\epsilon$-simply periodic curve $\zeta$ such that the two
endpoints of $\zeta$ are homoclinically related to $p.$

\end{lem}

\noindent{\bf Proof of Theorem B.} Let $\mathcal{U}(f)$ be a $C^1$
neighborhood  of $f\in{\mathcal{RT}(\mathbb{T}^3})$
 and let $f\in\mathcal{G}=\mathcal{G}_1\cap\mathcal{G}_2$.   Let $p$ be a hyperbolic periodic point of $f$.
 Then for any $g\in\mathcal{G}\cap\mathcal{U}(f)$, we have $\mathbb{T}^3=H_g(p_g)$, where $p_g$ is the continuation of $p.$
 Since $H_g(p_g)$ is not hyperbolic, from
\cite[section 4]{SV}, for any $\epsilon>0$, there is
$g_1\in\mathcal{G}\cap\mathcal{U}(f)$ such that $g_1$ has an
$\epsilon$-simply periodic curve $\eta$, whose endpoints are
homoclinically related to $p_{g_1}.$ 
Note that $\eta$ is a closed and connected set, so it is a
nontrivial continuum. Let $e=2\epsilon$ be an expansive constant
of $g_1.$ For $x\in \mathbb{T}^3$, we set
$\Gamma_e(x)=\{y\in\mathbb{T}^3:d(g_1^i(x), g_1^i(y))\leq e$ for
all $i\in\Z\}$. Since $g_1^{\pi(p_g)}(\eta)=\eta$, we have
$$L(g_1^{i\pi(p_g)}(\eta))=L(\eta)<e.$$

Clearly $\eta\subset\Gamma_e(x).$  Since $\eta$ is compact and
connected, and so, $\eta$ is not singleton which is
a contradiction.\hfill$\square$\\

\section*{Acknowledgement}

This work is supported by Basic Science Research Program through
the National Research Foundation of Korea (NRF) funded by the
Ministry of Science, ICT \& Future Planning
(No-2014R1A1A1A05002124).

\end{document}